\documentclass[12pt, a4paper]{amsart}
\usepackage{amssymb, amsmath, amsfonts, verbatim}
\newtheorem{lemma}[equation]{Lemma}
\newtheorem{thm}[equation]{Theorem}

\title{RECENT PROGRESS IN THE SYMMETRIC GENERATION OF GROUPS}
\author{Ben Fairbairn}
\address{Ben Fairbairn, Departmento de Mathematicas, Universidad de los Andes, Carrera 1 No 18A-12, Bogot\'{a}, Colombia}
\email{bt.fairbairn20@uniandes.edu.co}

\begin{document}

\maketitle

\begin{abstract}Many groups possess highly symmetric generating sets that are naturally endowed with an underlying combinatorial structure. Such generating sets can prove to be extremely useful both theoretically in providing new existence proofs for groups and practically by providing succinct means of representing group elements. We give a survey of results obtained in the study of these symmetric generating sets. In keeping with earlier surveys on this matter, we emphasize the sporadic simple groups.\end{abstract}

ADDENDUM: This is an updated version of a survey article originally accepted for inclusion in the proceedings of the 2009 `Groups St Andrews' conference \cite{orig}. Since \cite{orig} was accepted the author has become aware of other recent work in the subject that we incorporate to provide an updated version here (the most notable addition being the contents of Section \ref{Rubit}).

\section{Introduction}

This article is concerned with groups that are generated by highly symmetric subsets of their elements: that is to say by subsets of elements whose set normalizer within the group they generate acts on them by conjugation in a highly symmetric manner. Rather than investigate the behaviour of known groups we turn this procedure around and ask what groups can be generated by a set of elements that possesses a certain assigned set of symmetries. This enables constructions by hand of a number of interesting groups, including many of the sporadic simple groups. Much of the emphasis of the research project to date has been concerned with using these techniques to construct sporadic simple groups, and this article will emphasize this important special case. Recent work of the author and M\"{u}ller has been concerned with Coxeter groups, so we shall also describe this case too.

This article is intended as an `update' to the earlier survey article of Curtis \cite{ATLASconference}. Since \cite{ATLASconference} appeared several of the larger sporadic groups have succumbed to these techniques and a much wider class of reflections groups have been found to admit symmetric presentations corresponding to symmetric generating sets too. We refer the interested reader seeking further details to the recent book of Curtis \cite{Curtismain} which discusses the general theory of symmetric generation and many of the constructions mentioned here in much greater detail.

Throughout we shall use the standard {\sc Atlas} notation and conventions for finite groups and related concepts as defined in \cite{ATLAS}. 

We remark that the emphasis of this article is very much on symmetric generation in the sense of \cite{Curtismain}. Needless to say there also exist other slightly different, but closely related concepts of symmetric generation and progress in these areas have also taken place - for instance Abert's classification of symmetric presentations of abelian groups \cite{Abert}.

This article is organized as follows. In Section 2 we shall recall the basic definitions and notation associated with symmetric generation that are used throughout this article. In Section 3 we shall discuss symmetric generation of some of the sporadic simple groups. In Section 4 our attention will turn to reflection groups and in Section 5 we shall consider symmetric generators more general than cyclic groups of order 2, that earlier parts of the article ignore. After making some concluding remarks in Section 6 we finally give a short appendix describing some of the well-known properties of the Mathieu group M$_{24}$ that are used in several earlier sections.

\textbf{Acknowledgments} The author wishes to express his deepest gratitude to Dr John Bradley for bringing to his attention the recent work on the Rudvalis group that briefly outline in Section \ref{Rubit}. The author is also grateful to the referee of the original version of the article \cite{orig} for making several helpful comments and suggestions that have been incorporated here.

\section{Symmetric Generation}

\subsection{Definitions and Notation}
In this section we shall describe the general theory of symmetric generation giving many of the basic definitions we shall use throughout this survey.

Let $2^{\star n}$ denote the free group generated by $n$ involutions. We write $\lbrace t_1,t_2,\dots,t_n\rbrace$ for a set of generators of this free product. A permutation $\pi\in S_n$ induces an automorphism of this free product, $\hat\pi$, by permuting its generators namely \begin{equation}
t_i^{\hat{\pi}}:=\pi^{-1}t_i\pi=t_{\pi(i)}.\label{diamond}
\end{equation}
Given a group $N\leq S_n$ we can use this action to form a semi-direct product $\mathcal{P}$$:=2^{\star n} \colon N$.

When $N$ acts transitively we call $\mathcal{P}$ a \emph{progenitor}. Note that some of the early papers on symmetric generation insisted that $N$ acts at least 2-transitively. We do not impose this restriction here. Elements of $\mathcal{P}$ can all be written as a relator of the form $\pi w$ where $\pi\in N$ and $w$ is a word in the symmetric generators using (\ref{diamond}). Consequently any finitely generated subgroup of $\mathcal{P}$ may be expressed as $H:=\langle w_1\pi_1,\ldots,w_r\pi_r\rangle$ for some $r$. We shall express the quotient of $\mathcal{P}$ by the normal closure of $H$, $H^{\mathcal{P}}$, as
\begin{equation}
\frac{2^{\star n}\colon N}{w_1\pi_1,\ldots,w_r\pi_r}:=G.\label{spade}
\end{equation}
We say the progenitor $\mathcal{P}$ is factored by the relations $w_1\pi_1,\ldots,w_r\pi_r$. Whenever we write a relator $w\pi$ we shall tacitly be referring to the relation $w\pi=id$ thus we shall henceforth only refer to relations, when \emph{senso stricto} we mean relators. We call $G$ the \emph{target group}. Often these relations can be written in a more compact form by simply writing $(\pi w)^d$ for some positive integer $d$. It is the opinion of the author that no confusion should arise from calling both $t\in\mathcal{P}$ and its image in $G$ a \emph{symmetric generator}. Similarly no confusion should arise from calling both $N\leq\mathcal{P}$ and its image in $G$ the \emph{control group}. We define the \emph{length} of the relation $\pi w$ to be the number of symmetric generators in $w$.

Henceforth we shall slightly abuse notation in writing $t_i$ both for a generator of $2^{\star n}$ in $\mathcal{P}$ and for its homomorphic image in $G$. Similarly we shall write $N$ both for the control group in $\mathcal{P}$ and for the homomorphic image of $N$ in $G$. We shall also assume that $N$ is isomorphic to its image in $G$ as this is often the case in the most interesting examples. Again, it is the opinion of the author that no confusion should arise from this.

We are immediately confronted with the question of how to decide if $G$ is finite or not. To do this we resort to an enumeration of the cosets of $N$ in $G$. Let $g\in G$. We have that $gN\subset NgN$. Consequently the number of double cosets of the form $NgN$ in $G$ will be at most the number of single cosets of the form $gN$ in $G$ making them much easier to enumerate. To do this we make the following definition. 

Given a word in the symmetric generators, $w$, we define the \emph{coset stabilizing subgroup} to be the subgroup defined by\[ N^{(w)}:=\lbrace\pi\in N | Nw\pi=Nw\rbrace. \] This is clearly a subgroup of $N$ and there are $|N:N^{(w)}|$ right cosets of $N^{(w)}$ in the double coset $NwN$.

\subsection{Double Coset Enumeration}

As noted above, to verify a symmetric presentation of a finite group it is usual to perform a double coset enumeration. In early examples of symmetric presentations, the groups involved were sufficiently small for the coset enumeration to be easily performed by hand. However, attention has more recently turned to larger groups and consequently automation of the procedure to enumerate double cosets has been necessary.

Bray and Curtis have produced a double coset enumeration program specially suited to this situation in the {\sc Magma} computer package \cite{MAGMA}; it is described in \cite{cosetenumerator}. The program uses an adaptation of the celebrated Todd-Coxeter algorithm first described in \cite{ToddCoxeter} and follows an adaptation of an earlier program written by Sayed described in \cite[Chapter 4]{Sayid} which worked well with relatively small groups, but could not be made to cope with groups of a larger index or rank \cite[p.66]{Curtismain}.

\subsection{Some General Lemmata}
Given a particular progenitor we are immediately confronted with the problem of deciding what relations to factor it by. The following lemmata naturally lead us to relations that are often of great interest. In particular these lemmata prove surprisingly effective in naturally leading us to relations to consider.

The following lemma, given a pair of symmetric generators $\{t_i,t_j\}$, tells us which elements of the control group may be expressed as a word in the symmetric generators $t_i$ and $t_j$ inside the target group.

\begin{lemma}
If $G$ is a true image of the progenitor $\mathcal{P}= 2^{\star n} : N$ with   $\phi:P\rightarrow G$, and if $w(t_i, t_j)$ is in the kernel of  $\phi$ then $\pi$ must lie in the centralizer in $N$ of the stabilizer in 
$N$ of $i$ and $j$. That is
$$\pi\in C_N(N_{ij})$$,
where $N_{ij}$ denotes the stabilizer in $N$ of $i$ and $j$.
\end{lemma}

This lemma usual stated in the following more succinct fashion.

\begin{lemma}\label{famous}
$$\langle t_i, t_j \rangle\cap N\leq C_N(Stab_N(i,j))$$.
\end{lemma}

This lemma can easily be extended to an arbitrary number of symmetric generators by an obvious induction. Despite its strikingly minimal nature, this lemma proves to be extraordinarily powerful. For a proof see \cite[p.58]{Curtismain}.

Whilst the above lemma tells which elements of the control group may appear in a relation to factor a progenitor by, the precise length of this word remains open. A lemma that proves to be useful in settling this matter is the following.

\begin{lemma}\label{second}
Let $\mathcal{P}$$:=2^{\star n}:N$ be a progenitor in which the control group $N$ is perfect. Then any homomorphic image of $\mathcal{P}$ is either perfect or possesses a perfect subgroup to index 2. If $w$ is a word in the symmetric generators of odd length, then the image

$$\frac{2^{\star n}:N}{\pi w}$$
is perfect.
\end{lemma}

See \cite[Section 3.7]{Curtismain} for details. Other general results of this nature may also be found in \cite[Chapter 3]{Curtismain}.

\section{The Sporadic Simple Groups}

In this section we shall see how the techniques described in the previous section may be used to find symmetric presentations for some of the sporadic simple groups. We shall also see some examples in which symmetric presentations lead to new existence proofs for several of these groups as well as computational applications of providing the means of succinctly representing group elements. We proceed in ascending order of order.

\subsection{The Janko Group J$_3$}

In \cite{BradleyThesis} Bradley considered the primitive action of the group L$_2$(16):4 on 120 points. Since this action is transitive we can form the progenitor $2^{\star120}:(\mbox{L}_2(16):4)$. Using a single short relation found by Bray, Bradley was able to verify the symmetric presentation$$\frac{2^{\star120}:(\mbox{L}_2(16):4)}{(\pi t_1)^5}\cong\mbox{J}_3:2$$where $\pi\in$ L$_2$(16):4 is a well chosen permutation of order 12 by performing by hand an enumeration of the double cosets of the form (L$_2(16):4)w(\mbox{L}_2(16):4)$ where $w$ is a word in the symmetric generators.

The full double coset enumeration in \cite{BradleyThesis} is a fairly long and involved calculation. In \cite{J3} Bradley and Curtis state and prove the general lemmata used to perform this calculation and describe the Cayley graph for the symmetric generating set derived from it. They go on to use this symmetric presentation to provide a new existence proof for J$_3$:2 by proving that the target group in this symmetric presentation must either have order 2 or the correct order to be J$_3$:2. They then exhibit $9\times9$ matrices over the field of four elements that satisfy the relations of the presentation, thus verifying the above isomorphism.

An immediate consequence of this presentation is that elements of J$_3$:2 may be represented as an element of L$_2$(16):4 followed by a short word in the symmetric generators. In \cite{BradleyThesis} Bradley gives a program in the algebra package {\sc Magma} \cite{MAGMA} for multiplying elements represented in this form together and expressing their product in this concise form.

\subsection{The McLaughlin Group McL}\label{Mcluff}

The Mathieu group M$_{24}$ naturally acts on the set of 2576 dodecads (see Appendix). If we fix two of the 24 points that M$_{24}$ acts on then there are 672 dodecads containing one of these points but not the other. The stabilizer of the two points, the Mathieu group M$_{22}$, acts transitively on these 672 dodecads. We can therefore form the progenitor $2^{\star672}:\mbox{M}_{22}$. In \cite{BradleyThesis,McL} Bradley anbd Curtis investigated this progenitor and was naturally led to the symmetric presentation$$\frac{2^{\star672}:\mbox{M}_{22}}{\pi(t_At_B)^2}\cong\mbox{McL:2}$$where $A$ and $B$ are two dodecads intersecting in eight points and $\pi\in\mbox{M}_{22}$ is a permutation arrived at using Lemma\ref{famous}. As with J$_3$:2 the coset enumeration in this case was performed entirely by hand, thus providing a new computer-free existence proof for this group.

In \cite{McL} Bradley and Curtis go on to consider the subprogenitor $2^{\star42}:\mbox{A}_7$ defined by the natural action of A$_7$ on the seven points of a fixed heptad and on a certain orbit of dodecads in the complement. Motivation for an additional relation needed to define the target group comes from a symmetric presentation for the unitary group U$_3(5):2$ described in terms of Hoffman-Singleton graph, the details of which we ommit. We are thus lead to the following symmetric presentation.
$$\frac{2^{\star42:\mbox{A}_7}}{(t_st_t)^2=(36)(45)}\cong\mbox{M}_{22}$$

\subsection{The Conway Group $\cdot$0=2\.{}Co$_1$}\label{dotto}In \cite{revisited} Bray and Curtis consider the progenitor $2^{\star{24 \choose 4}}:\mbox{M}_{24}$ defined by natural the action of the Mathieu group M$_{24}$ on subsets of size four of a set of order 24 on which M$_{24}$ naturally acts (see Appendix). After eliminating words of length 2 and other relations of length three they are immediately led to the symmetric presentation
\[\frac{2^{\star{24\choose4}}\colon\mbox{M}_{24}}{\pi t_{ab}t_{ac}t_{ad}}\cong\cdot0,\]
where $a,b,c$ and $d$ are pairs of points the union of which is a block of the $\mathcal{S}$(5,8,24) Steiner system on which M$_{24}$ naturally acts (see the {\sc Atlas}, \cite[p.94]{ATLAS} and the Appendix) and $\pi\in$ M$_{24}$ is the unique non-trivial permutation of M$_{24}$ determined by the Lemma \ref{famous}. From this presentation, an irreducible 24 dimensional $\mathbb{Z}$ representation is easily found and considering the action this gives on certain vectors in $\mathbb{Z}^{24}$ the famous Leech lattice effortlessly `drops out'. Furthermore the symmetric generators are revealed to be essentially the elements of $\cdot$0 discovered by Conway when he first investigated the group \cite{ConOrig}.

Using this symmetric presentation of $\cdot$0 and detailed knowledge of the coset enumeration needed to verify it Curtis and the author have been able to produce a program in {\sc Magma} \cite{MAGMA}, available from the author's website (which at the timeof writing may be found at

\begin{center}
\texttt{http://matematicas.uniandes.edu.co/}$\sim$\texttt{benfairbairn/Homepage.htm}), 
\end{center}

that represents elements of$\cdot$0 as a string of at most 64 symbols and typically far fewer, essentially by representing them as $\pi w$ where $\pi\in\mbox{M}_{24}$ and $w$ is a word in the symmetric generators \cite{CurtisFairbairn}. This represents a considerable saving compared to representing an element of $\cdot$0 as a permutation of196560 symbols or as a $24 \times 24$ matrix (ie as a string of 24$^2$=576 symbols). To date, $\cdot$0 is the largest group for which a program of this kind has been produced.

\subsection{The Rudvalis Group Ru}\label{Rubit}

In \cite{Ru} Bradley, Curtis and Aslam Malik exhibit two symmetric presentations of the Rudvalis group, rectifying a notable omission from \cite{Curtismain}.

For the first symmetric presentation we consider an imprimitive degree 105 action of the linear group L$_4$(2). The group L$_4$(2) acts naturally on the 15 non-zero vectors of a 4-dimensional $\mathbb{F}_2$ vector space $V$. This induces a primitive action on the ${15\choose3}/3=35$ 2-dimensional subspaces. Since the stabilizer of a 2 dimensional subspace $W\leq V$ is transitive on the three 1-dimensional subspaces that are contained in $W$, we have a transitive action on 105 points. The authors of \cite{Ru} call these configurations \emph{matchsticks} and using the action of L$_4$(2) on the set of all matchsticks we define the progenitor $2^{\star105}:\mbox{L}_4(2)$ leading to the symmetric presentation
$$\frac{2^{\star105}:\mbox{L}_4(2)}{t_At_Et_C=1,t_At_Bt_Ct_D=\sigma}\cong\mbox{Ru}$$
where the relationships between the matchsticks $A,\ldots,E$ is too complicated to describe here and the element $\sigma\in\mbox{L}_4(2)$ being an element that is naturally arrived at using Lemma \ref{famous}. We remark that the above symmetric presentation is not explicitly written down in \cite{Ru}, but it is shown to be equivalent to a symmetric presentation defined using the same progenitor that, whilst slightly easier to describe, requires three relations instead of two and is thus slightly more difficult to naturally motivate - see \cite[p.11]{Ru}.

For the second symmetric presentation of Ru our attention turns to the natural action of the linear group L$_3$(2) on 7 points and the progenitor $2^{\star7}:\mbox{L}_3(2)$ defined by this action. This progenitor has proved useful in providing interesting symmetric presentations before - the Mathieu group M$_{24}$ being constructed as a homomorphic image of this progenitor - see \cite{Mathieu} (\cite[Section I]{Curtismain} also contains a detailed discussion of this presentation, the cover picture of \cite{Curtismain} being a diagrammatic illustration showing this presentation's relationship to the celebrated Klein quartic curve). 

The motivation for this presentation comes from an earlier symmetric presentation for the Tits group $^2\mbox{F}_4(2)$ that was first met in passing during the `systematic approach' of Curtis, Hammas and Bray in \cite{syst}. The natural action of the symmetric group on four points enables us to define the progenitor $2^{\star4}:\mbox{S}_4$ leading to the symmetric presentation 
$$\frac{2^{\star4}:\mbox{S}_4}{((1234)t_1)^{10},((123)t_1)^{13},((1234)t_3t_2[t_1,t_2t_3])^3}\cong\mbox{ }^2\mbox{F}_4(2).$$
Treating the progenitor $2^{\star4}:\mbox{S}_4$ as a subprogenitor of $2^{\star7}:\mbox{L}_3(2)$ the above symmetric presentation can be extended to the following symmetric presentation of the full Rudvalis group.
$$\frac{2^{\star7}:\mbox{L}_3(2)}{((14)(0563t_0)^{10}, ((146)(253)t_0)^{13},((14)(0563)t_0t_3[t_6,t_3t_0])^3, }\cong\mbox{Ru}$$
$$((14)(0563)t_1(t_0t_6)^2t_1)^2, ((14)(0563)t_2(t_1t_4)^2)^2$$

\subsection{The Janko Group J$_4$}\label{Janko}

In \cite{J4paper} Bolt, Bray and Curtis, building on earlier work of Bolt \cite{BoltThesis}, considered the primitive action of the Mathieu group M$_{24}$ on the 3795 triads (see Appendix). Since this action is transitive we can define the progenitor $2^{\star3795}:\mbox{M}_{24}$ leading to the presentation
$$\frac{2^{\star3795}:\mbox{M}_{24}}{t_At_Bt_C,\pi t_At_Dt_At_E}\cong\mbox{J}_4$$
where the triads $A$, $B$ and $C$ each have a single octad in common, the remaining octads intersecting in either four or no points, depending on whether they belong to the same triad or not. The triads $A$, $D$ and $E$ have a more complicated relationship to one another and we merely refer the reader to \cite[p.690]{J4paper} for details. Again, the permutation $\pi\in\mbox{M}_{24}$ is naturally arrived at by Lemma \ref{famous}.

Note that in this case the `excess' of relations is merely an illusion - removing the first of the two relations gives a symmetric presentation of the group J$_4\times2$, a relation of odd length being needed to produce a simple group in accordance with Lemma \ref{second}.

\subsection{The Fischer Groups}

In the earlier survey of Curtis \cite{ATLASconference} it was noted that the each of the sporadic Fischer groups were homomorphic images of progenitors defined using non-involutory symmetric generators and that relations defining the Fischer groups were at the time of writing being investigated by Bray.

Since that time involutory symmetric presentations for each of the sporadic Fischer groups as well as several of the classical Fischer groups closely related to them have been found. In \cite{Fischer1} Bray, Curtis, Parker and Wiedorn proved
$$\frac{2^{\star{10\choose4}}:\mbox{S}_{10}}{((45)t_{1234})^3,(12)(34)(56)t_{1234}t_{1256}t_{3456}t_{7890}}\cong\mbox{Sp}_8(2),$$
$$\frac{2^{\star(7+{7\choose3})}:\mbox{S}_7}{((12)t_1)^3,((45)t_{1234})^3,(12)(34)(56)t_{1234}t_{3456}t_{1256}t_7}\cong\mbox{Sp}_6(2),$$
$$\frac{2^{\star288}:\mbox{Sp}_6(2)}{(r_{4567}t_{\infty})^3}\cong3\mbox{\.{}O}_7(3).$$
In the first case the progenitor is defined by the natural action of the symmetric group S$_{10}$ on subsets of size four.  The second `progenitor', which deviates from the traditional definition of progenitor since the control group does not act transitively, is defined by the natural action of the symmetric group S$_7$ on seven points and its action on subsets of size three. In the third case the action defining the progenitor is the action of Sp$_6$(2) on the cosets of a subgroup isomorphic to the symmetric group S$_7$. The symbol $r_{4567}$ appearing in the third of these presentations corresponds to a symmetric generator defined in the above symmetric presentation of Sp$_8$(2).

In \cite{Fischer2} the same authors, motivated by the above results were able to exhibit symmetric presentations for the sporadic Fischer groups as follows.

For a symmetric presentation of the sporadic Fischer group Fi$_{22}$ we consider the action of the group $2^6:\mbox{Sp}_6(2)$ (recalling the well known symmetric group isomorphism Sp$_6(2)\cong\mbox{S}_6$) and using its action on the cosets of a copy of a subgroup isomorphic to the symmetric group S$_8$ which leads us to the presentation

$$\frac{2^{\star2304}:(2^6:\mbox{Sp}_6(2))}{(ts)^3}\cong\mbox{3\.{}Fi}_{22}$$

\noindent where $t$ is a symmetric generator and $s$ is a well chosen element of the control group the choice of which is motivated by the symmetric presentations for the classical Fischer groups given above.  

For a symmetric presentation of the sporadic Fischer group Fi$_{23}$ we consider the natural action of the symmetric group S$_{12}$ on the partitions into three subsets size four of the set $\{1,\ldots,9,0,x,y\}$ . Using this we can define the progenitor $2^{\star5775}:\mbox{S}_{12}$ which we factor by the relations

\begin{center}
$\left(\begin{array}{c|c|}
&\underline{1234}\\
(15)&\underline{5678}\\
&90xy\\ \end{array}\hspace{1mm}\right)^3$
\end{center}

\noindent and 

\begin{center}
$\begin{array}{c|c||c||c|}
&\underline{1234}&\underline{1256}&\underline{1278}\\
(12)(34)(56)(78)&\underline{5678}&\underline{3478}&\underline{3456}\\
&90xy&90xy&90xy\\ 
\end{array}.$
\end{center}

The above presentation was set-up to resemble the classical relations satisfied by the bifid maps related to the Weyl groups of type E$_6$ and E$_7$ (see \cite[Section 4.4]{Curtismain} for a discussion relating these to symmetric presentations). In addition to the above presentation, the same authors go on to prove another, substantially simpler, symmetric presentation for Fi$_{23}$ namely 
$$\frac{2^{\star13056}:\mbox{Sp}_8(2)}{(r_{1234}t_1)^3}\cong\mbox{Fi}_{23}$$
where the action defining the progenitor in this case is the action of the symplectic group Sp$_8$(2) on the cosets of the maximal subgroup isomorphic to the symmetric group S$_{10}$ and $r_{1234}$ is the symmetric generator defined in the symmetric presentation of Sp$_8$(2) given earlier.

Finally for a symmetric presentation of the sporadic Fischer group Fi$_{24}$ we consider the action of the orthogonal group O$_{10}^-(2):2$ on the cosets of a copy of a subgroup isomorphic to the symmetric group S$_{12}\leq\mbox{O}_{10}^-:2$ which leads us to the symmetric presentation

 $$\frac{2^{\star104448}:(\mbox{O}_{10}^-(2):2)}{(ts)^3}\cong\mbox{3\.{}Fi}_{24}$$
 
\noindent where $t$ is a symmetric generator and $s$ is a well chosen element of the control group the choice of which is motivated by the symmetric presentations for the classical Fischer groups given above.

\section{Coxeter Groups}

Recall that a \emph{Coxeter diagram} of a presentation is a graph in which the vertices correspond to involutory generators and an edge is labeled with the order of the product of its two endpoints.Commuting vertices are not joined and an edge is left unlabeled if the corresponding product has order three. A Coxeter diagram and its associated group are said to be \emph{simply laced} if all the edges of the graph are unlabeled. In \cite{ATLASconference} Curtis notes that if such a graph has a ``tail'' of length at least two, as in Figure I, then we see that the generator corresponding to the terminal vertex, $a_r$, commutes with the subgroup generated by the subgraph $\mathcal{G}_0$.

\newpage
\setlength{\unitlength}{1mm}
\begin{picture}(0,0)(-15,0)
\put(30,0){\circle{20}}\put(37,0){\line(1,0){40}}\multiput(57,0)(20,0){2}{\circle*{2}}
\put(28.5,-1){$\mathcal{G}_0$}\put(56,-5){$a_{r-1}$}\put(76,-5){$a_r$}\end{picture}
\bigskip\bigskip
\begin{center}Figure I: A Coxeter diagram with a tail.
\end{center}

\bigskip

The author has more recently investigated a slight generalization of this idea to produce extremely succinct symmetric presentations for all of the finite simply laced Coxeter groups. More specifically in \cite[Chapter 3]{nearly, thesis} we prove

\begin{thm} \label{main} Let $S_n$ be the symmetric group acting on $n$ objects and $W(\Phi)$ denote the Weyl group of the root system $\Phi$. Then:

\begin{enumerate}
\item\[ \frac{2^{\star {n \choose 1}}:S_n}{(t_1(12))^3}\cong W(\mbox{A}_{n}) \]
\item\[ \frac{2^{\star {n \choose 2}}:S_n}{(t_{12}(23))^3}\cong W(\mbox{D}_n) \mbox{ for } n\geq4\]
\item\[ \frac{2^{\star {n \choose 3}}:S_n}{(t_{123}(34))^3}\cong W(\mbox{E}_{n}) \mbox{ for }n=6,7,8.\]
\end{enumerate}
\end{thm}

In particular, we show that the above symmetric presentations may be naturally arrived at using general results such as Lemma \ref{famous} without considering the general theory of Coxeter groups. Subsequently, the author and M\"{u}ller generalized this result as follows. Let $\Pi$ be a set of fundamental reflections generating a Coxeter group $G$. If $s\in\Pi$ is a reflection, the subgroup $N\leq G$ generated by $\Pi\backslash\{s\}$ may be used as a control group for $G$ and the action of $N$ on the set $\{s^N\}$ may be used to define a progenitor, $\mathcal{P}$. Factoring $\mathcal{P}$ by the relations defining $G$ that involve $s$ provides a symmetric presentation of $G$. The lack of restrictions on $G$ means it is possible that $G$ is infinite. In particular our theorem applies to affine reflection groups and to many of the hyperbolic reflection groups. Our theorem thus provides the first examples of interesting infinite groups for which a symmetric presentation has been found (`interesting' in the sense that any progenitor is an infinite group that is symmetrically generated \emph{ab initio}, but only in unenlightening manner). See \cite{muller} for details.

We remark that in \cite{muller} various (extremely weak) finiteness assumptions about the control group are made to ensure that the progenitors considered there do actually satisfy the definition of progenitor. The proofs of the results given there, however, nowhere require the progenitor to be finitely generated or the action defining the progenitor to be on a finite number of points. If the definition of progenitor is weakened slightly to allow this wider class of objects to be considered then the results given in \cite{muller} and their proofs remain valid.

\section{Non-Involutory Symmetric Generation}

Whilst many of the applications of symmetric generation have focused on the case where the symmetric generators are involutions, it is also possible to consider groups in which the symmetric generating set contains groups more general than cyclic groups of order 2.

Let $H$ be a group. Instead of considering the free group $2^{\star n}$ we consider the free product of $n$ copies of $H$, which we shall denote $H^{\star n}$. If $N$ is a subgroup of $S_n\times Aut(H)$ then $N$ can act on $H^{\star n}$ and we can form a progenitor in a similar manner to before usually writing $H^{\star n}:_m N$ to indicate that the action defining the progenitor is not necessarily just a permutation of the copies of $H$ - the `$m$' standing for `monomial' - see below. The special case in which $H$ is cyclic of order 2 is simply the involutory case discussed in earlier sections.

One common source of actions used for defining progenitors in this way is as follows. Recall that a matrix is said to be \emph{monomial} if it has only one non-zero entry in each row and column. A representation of a group $\rho:G\rightarrow GL(V)$ is said to be \emph{monomial}, if $\rho(g)$ is a monomial matrix for every $g\in G$.

An $n$ dimensional monomial representation of a control group $N$ may be used to define a progenitor $H^{\star n}:_m N$, each non-zero entry of the monomial matrix acting as a cyclic subgroup of $Aut(H)$. 

In \cite{WhyteThesis} Whyte considered monomial representations of decorations of simple groups and their use in forming monomial progenitors. In particular Whyte classified all irreducible monomial representations of the symmetric and alternating groups (reproducing earlier work of Djokovi\'{c} and Malzan \cite{Malzan, Malzan2}) and their covers. Whyte went on to classify the irreducible monomial representations of the sporadic simple groups and their decorations and found a large number of symmetric presentations defined using the actions defined by these monomial representations on free products of cyclic groups \cite{spormono}.

In light of Lemma \ref{second}, it is natural to restrict our attention to control groups that are perfect, and in particular control groups that are simple. In light of the classification of the finite simple groups, the work of Whyte only left open the many families of groups of Lie type. Since the most `dense' family of groups of Lie type are the groups $L_2(q)$ it is natural to consider monomial representations of these groups. In \cite{linearpap} using only elementary techniques the author classifies the monomial representations for all the natural decorations of $L_2(q)$ and these where put to great use in defining symmetric presentations in \cite[Chapter 4]{thesis}.

\section{Concluding Remarks}

We note that the symmetric presentations discussed so far represent the most `elegant' symmetric presentations of sporadic simple groups discovered since \cite{ATLASconference}. Several other symmetric presentations of sporadic groups that have proved more difficult to motivate for various reasons have also been discovered. Examples include Bray's presentation of the Lyons group Ly \cite[Chapter10]{BrayThesis}, Bolt's presentation of the Conway group Co$_3$ \cite[Section 3.3]{BoltThesis} and the author's presentation of the Conway group Co$_2$ (unpublished). Furthermore several of the larger sporadic groups, most notably the Thompson group Th, the Baby Monster $\mathbb{B}$ and of course the Monster $\mathbb{M}$, still lack any kind of symmetric presentation at all, though various conjectures in these cases do exist (see \cite[Section 4.1]{ATLASconference}).

\section{Appendix: The Mathieu Group M$_{24}$}

In this appendix we gather together some well-known facts about the Mathieu group M$_{24}$ that are repeatedly called on in Section 3. Let $X$ be a set such that $|X|=24$. Recall that the Steiner system $\mathcal{S}$(5,8,24) is a certain collection of 759 subsets $O\subset X$ with $|O|=8$ called \emph{octads} such that any $P\subset X$ with $|P|=5$ is contained in a unique octad. An automorphism of $\mathcal{S}$(5,8,24) is a permutation of $X$ such that the set of all octads is preserved. The group of all automorphisms of$\mathcal{S}$(5,8,24) is the sporadic simple Mathieu group M$_{24}$. The group M$_{24}$ and the Steiner system $\mathcal{S}$(5,8,24) have many interesting properties, some of which we summarize below.

\begin{itemize}\item The pointwise stabilizer of two points of $X$ is another sporadic simple group known as the Mathieu group M$_{22}$ - see Section \ref{Mcluff}.
\item If $A$ and $B$ are two octads such that $|A\cap B|=2$ then their symmetric difference $(A\cup B)\backslash(A\cap B)$ is called a \emph{dodecad}. There are 2576 dodecads and M$_{24}$ acts transitively on the set of all dodecads - see Section \ref{Mcluff}.
\item The group M$_{24}$ acts 5-transitively on the 24 points of $X$ (so in particular acts transitively on subsets of $X$ of cardinality 4 -see Section \ref{dotto}).
\item The set $X$ may be partitioned into three disjoint octads. Such a partition is called a \emph{trio}. There are 3795 trios and M$_{24}$ acts primitively on the set of all trios - see Section \ref{Janko}.
\end{itemize}  

There are several good discussions of the Mathieu groups and their associated Steiner systems in the literature, many of which contain proofs of the above. See Conway and Sloane \cite[Chapters 10 \& 11]{SPLAG} or the {\sc Atlas} \cite[p.94]{ATLAS} and the references given therein for further details.

\end{document}